\documentclass{ifacconf}

\usepackage[dvips]{graphicx}
\usepackage{amssymb}
\usepackage{amsfonts}
\usepackage{amsmath}
\usepackage[mathscr]{eucal}
\usepackage{fancyhdr}
\usepackage{psfrag}    
\usepackage{datetime}       

\usepackage{natbib}

\renewcommand{\overline}{\bar}
\renewcommand{\limits}{}

\renewcommand{\L}{\mathcal{L}}
\renewcommand{\S}{\mathcal{S}}
\newcommand{\D}{\mathcal{D}}
\renewcommand{\d}{\text{d}}
\renewcommand{\phi}{\varphi}

\begin{document}
\begin{frontmatter}

\title{Memory Elements: A Paradigm Shift in Lagrangian Modeling of Electrical Circuits} 

\author{Dimitri Jeltsema}

\address{Delft Institute of Applied Mathematics, Delft University of Technology, Mekelweg 4, 2628 CD Delft, The Netherlands\\ (Email: d.jeltsema@tudelft.nl)}


\begin{abstract}
Meminductors and memcapacitors do not allow a Lagrangian formulation in the classical sense since these elements are nonconservative in nature and the associated energies are not state functions. To circumvent this problem, a different configuration space is considered that, instead of the usual loop charges, consist of time-integrated loop charges. As a result, the corresponding Euler-Lagrange equations provide a set of integrated Kirchhoff voltage laws in terms of the element fluxes. Memristive losses can be included via a second scalar function that has the dimension of action. A dual variational principle follows by considering variations of the integrated node fluxes and yields a set of integrated Kirchhoff current laws in terms of the element charges. Although time-integrated charge is a somewhat unusual quantity in circuit theory, it may be considered as the electrical analogue of a mechanical quantity called absement. Based on this analogy, simple mechanical devices are presented that can serve as didactic examples to explain memristive, meminductive, and memcapacitive behavior.\\[-1.5em]
\end{abstract}
\begin{keyword}
Memristor, Meminductor, Memcapacitor, Memory elements, Lagrangian modeling.\\[-0.5em]
\end{keyword}
\end{frontmatter}

\section{Introduction}

Memristors, meminductors, and memcapacitors constitutive an increasingly important class of two-terminal circuit elements whose resistance, inductance, and capacitance retain memory of the past states through which the elements have evolved. All three elements are nonlinear and can be identified by their pinched hysteresis loop in the voltage versus current plane, current versus flux plane, and voltage versus charge plane, respectively. While there are many discovered experimental realizations and applications of systems that exhibit memristive behavior, ranging from applications in non-volatile nano memory to intelligent machines with learning and adaptive capabilities, the number of systems showing memcapacitive and meminductive behavior is still somewhat limited. Nevertheless, several applications of these concepts are foreseen in the field of logic and arithmetic operations using memristive and memcapacitive devices, and field-programmable quantum computation using meminductive and memcapacitive devices. See \cite{DiVentra2009} for an introduction into the memory elements, and \cite{PershinMemoryEffects2011} for a fairly complete overview of the current state-of-the-art and an extensive list of references. 

In the past century, a significant amount of research has been devoted to the description and analysis of conventional electrical circuits using Lagrangian and Hamiltonian methods; see \cite{CSMpaper} and the references therein. Since these classical frameworks are important not just for their broad range of applications, but also for their role in advancing deep understanding of physics, the next step is to consider circuits made from memristors, meminductors, and memcapacitors. \newpage

~\\[0.01em]
The remainder of the paper is organized as follows. In Section \ref{sec:memory-defs}, the mathematical properties of the memristor, meminductor, and memcapacitor are discussed. The structure of the Lagrangian equations for a large class of conventional circuits is briefly reviewed in Section \ref{sec:EL-conventional} and it will be shown that meminductors and memcapacitors do not fit into this framework a priori since these elements are nonconservative in nature and the associated energies are not state functions. To circumvent this problem, in Section \ref{sec:EL-memory} a different configuration space is considered that, instead of the usual loop charges, consist of time-integrated loop charges. The Lagrangian is defined by the difference between two novel state functions in a fashion similar to the usual magnetic co-energy minus electric energy setup, but having the dimensions of energy times time-squared which, in turn, is equivalent to action times time. As a result, the corresponding Euler-Lagrange equations provide a set of integrated Kirchhoff voltage laws (iKVL's) in terms of the element fluxes. Memristive and resistive losses can be included via the introduction of a second scalar function that has the dimension of action. Furthermore, a dual variational principle follows by considering variations of the integrated node fluxes and yields a set of integrated Kirchhoff current laws (iKCL's) in terms of the element charges. Although time-integrated charge, which in SI units is measured in Coulomb times seconds, is a somewhat unusual quantity in circuit theory, it may be considered as the electrical analogue of a mechanical quantity called absement. Based on this analogy, in Section \ref{sec:absement}, simple mechanical devices are presented that can serve as didactic examples to explain meminductive and memcapacitive behavior. Finally, the paper is concluded with some final remarks in Section \ref{sec:conclusions}.

\section{Definition of Memory Elements}\label{sec:memory-defs}

From a mathematical perspective, the behavior of a two-terminal resistor, inductor, and capacitor, whether linear or nonlinear, is described by a relationship between two of the four basic electrical variables, namely, voltage $V$, current $I$, charge $q$, and flux $\phi$, where 
\begin{align}
   q(t) &= \int\limits_{-\infty}^{t} I(\tau) \d \tau,\\
\phi(t) &= \int\limits_{-\infty}^{t} V(\tau) \d \tau.
\end{align}
A resistor is described by a constitutive relationship between current and voltage; a capacitor by that of voltage and charge; and an inductor by that of current and flux linkage. 

Based on logical and symmetry reasonings, \cite{Chua1971} postulated the existence of a fourth element that is characterized by a constitutive relationship between the remaining two variables, namely charge and flux. This element was coined memristor (a contraction of memory and resistor) referring to a resistor with memory. The memory aspect stems from the fact that a memristor `remembers' the amount of current that has passed through it together with the total applied voltage. More specifically, if $q$ denotes the charge and $\phi$ denotes the flux, then a charge-modulated memristor is defined by the constitutive relationship $\phi=\hat{\phi}(q)$. Since flux is defined by the time integral of voltage $V$ (like in Faraday's law), and charge is the time integral of current $I$, or equivalently, $V=\dot{\phi}$ and $I=\dot{q}$, we obtain 
\begin{equation}\label{eq:e-Ohmmemris}
V=R_M(q)I,
\end{equation} 
where $R_M(q):=\d\phi/\d q$ is the incremental memristance. 

Note that (\ref{eq:e-Ohmmemris}) is the definition of a memristor in impedance form. The admittance form $I=G_M(\phi)V$, with incremental memductance $G_M(\phi):=\d q / \d \phi$, is obtained by starting from a constitutive relationship $q=\hat{q}(\phi)$.   

In addition to the memristor, it is shown in \cite{DiVentra2009} that the memory-effect can be associated to inductors and capacitors as well. For that, let $\sigma$ and $\rho$ denote the time-integrals of charge and flux, 
\begin{align}
\sigma(t) &= \int\limits_{-\infty}^{t} q(\tau) \text{d} \tau,\label{eq:sigma}\\
  \rho(t) &= \int\limits_{-\infty}^{t} \phi(\tau) \text{d} \tau,\label{eq:rho}
\end{align}
respectively. Then, a memory inductor, or meminductor for short, is a two-terminal element defined by a constitutive relation $\rho=\hat{\rho}(q)$. Indeed, differentiation of the latter with respect to time yields 
\begin{align}\label{eq:memL}
\phi=L_M(q)I.
\end{align}
Since $L_M(q):=\d\rho/\d q$ relates flux with current, its values clearly correspond to the units of inductance. Bearing in mind that charge is the time integral of current, the memory aspect of a charge-modulated meminductor stems from the fact that it `remembers' the amount of current that has passed through it.  

Dually, a flux-modulated memcapacitor is defined by a constitutive relationship $\sigma=\hat{\sigma}(\phi)$, which, after differentiation with respect to time, yields 
\begin{align}\label{eq:memC}
q=C_M(\phi)V,
\end{align}
where $C_M(\phi):=\d \sigma / \d \phi$ denotes the incremental capacitance. The memory aspect of a flux-modulated memcapacitor stems from the fact that it `remembers' the amount of voltage that has been applied to it. 

A meminductor (resp. memcapacitor) that depends on the history of its flux (resp. charge) can be formulated by starting from a constitutive relationship of the form $q=\hat{q}(\rho)$ (resp. $\phi=\hat{\phi}(\sigma)$). 

In the special case that the constitutive relationship of a memristor is linear, a memristor becomes an ordinary linear resistor. Indeed, in such case (\ref{eq:e-Ohmmemris}) reduces to $\phi=R q$, with constant memristance M (the slope of the line), or equivalently, $V=RI$, which precisely equals OhmÕs law. Hence it is not possible to distinguish a two-terminal linear memristor from a two-terminal linear resistor. The same holds for a linear meminductor and linear memcapacitor, where (\ref{eq:memL}) reduces to $\rho=Lq$, or equivalently, $\phi=LI$, and (\ref{eq:memC}) to $\sigma=C\phi$, or equivalently, $q=CV$, respectively.

\section{Self-Adjointness of Circuit Dynamics}\label{sec:EL-conventional}

The dynamical behavior of any electrical circuit consisting of conventional, possibly nonlinear, resistors, inductors, and capacitors is basically determined by three types of equations: those arising from Kirchhoff's voltage law (KVL), those arising from Kirchhoff's current law (KCL), and the constitutive relationships of the elements. In many cases this leads to differential equations of the form\footnote{Throughout the document we adopt the summation convention of repeated indices.}	
\begin{equation}\label{eq:general-circuit}
A_{ij}(\dot{x})\ddot{x}^j + B_i(x,\dot{x})=0,   
\end{equation}
where $i,j=1,\ldots,n$, and $x \in \mathbb{R}^n$ represents a column vector of loop charges and node fluxes. 

The system of differential equations (\ref{eq:general-circuit}) allow a Lagrangian description if we can find a Lagrangian $\L(x,\dot{x})$ satisfying
\begin{equation}\label{eq:IP}
\frac{\d}{\d t}\frac{\partial \L}{\partial \dot{x}^i} - \frac{\partial \L}{\partial x^i} \equiv A_{ij}(\dot{x})\ddot{x}^j + B_i(x,\dot{x}).
\end{equation} 
In mechanics it is known that the existence of a Lagrangian $\L(x,\dot{x})$ relies on the fact that the system of differential equations (\ref{eq:general-circuit}) is self-adjoint, which for the present form is tantamount to the following set of integrability conditions \citep{Santilli1978}:
\begin{align}
A_{ij} &= A_{ji},\label{eq:SA1}\\[0.66em]
\frac{\partial A_{ik}}{\partial \dot{x}^j} &= \frac{\partial A_{jk}}{\partial \dot{x}^i},\label{eq:SA2}\\[0.3em]
\frac{\partial B_i}{\partial x^j} - \frac{\partial B_j}{\partial x^i} &= \frac{1}{2} \frac{\partial}{\partial x^k} {\left(\frac{\partial B_i}{\partial \dot{x}^j} - \frac{\partial B_j}{\partial \dot{x}^i}  \right)}\dot{x}^k ,\label{eq:SA3}\\[0.3em]
\frac{\partial B_i}{\partial \dot{x}^j} + \frac{\partial B_j}{\partial \dot{x}^i} &= 0.\label{eq:SA4}
\end{align}

\subsection{Conventional Conservative Circuits}\label{subsec:EL-LC}

From a Lagrangian perspective, it is well-known that for a large class of  circuits that are made from current-controlled inductors and charge-controlled capacitors, and consisting of $n$ loops, the differential equations (\ref{eq:general-circuit}) can be written as
\begin{equation}\label{eq:EL-conservative}
\frac{\d}{\d t}\frac{\partial \L}{\partial \dot{q}^i} - \frac{\partial \L}{\partial q^i} = 0,
\end{equation}
where $q^i$ represents the loop charge associated to the $i$-th loop, and $\dot{q}^i$ represents the loop current circulating in the $i$-th loop. For this particular case, the Lagrangian $\L(q,\dot{q})$ equals the total magnetic co-energy stored in the inductors, $T^*(\dot{q})$, minus the total electric energy stored in the capacitors, $U(q)$.

On the other hand, if the inductors are flux-controlled and the capacitors are voltage-controlled, we need to start from a node analysis yielding a so-called co-Lagrange equation associated to each node in the network. Hence, if $\dot\phi^j$ represents the potential of the $j$-th node together with its time-integral, the node flux $\phi^j$, one obtains
\begin{equation}\label{eq:EL*-conservative}
\frac{\d}{\d t}\frac{\partial \L^*}{\partial \dot{\phi}^j} - \frac{\partial \L^*}{\partial \phi^j} = 0,
\end{equation}
where the co-Lagrangian $\L^*(\phi,\dot\phi)$ equals the total electric co-energy stored in the capacitors, $U^*(\dot{\phi})$, minus the total magnetic energy stored in the inductors, $T(\phi)$.

\begin{figure}[h]
\begin{center}
\psfrag{q}[][]{$q$}
\psfrag{p}[][]{$\phi$}
\psfrag{i}[][]{$I$}
\psfrag{C}[][]{$C$}
\psfrag{L}[r][]{$L$}
\includegraphics[height=30mm]{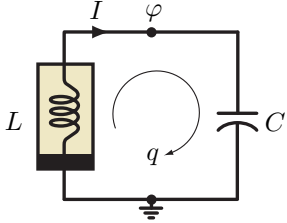}
\caption{Circuit with a nonlinear conventional inductor.}
\label{fig:LCexample}
\end{center}
\end{figure}

As an illustration, consider a circuit consisting of a nonlinear current-controlled inductor, with constitutive relation $\phi=\hat\phi(I)$, and a linear capacitor $C$ as shown in Fig.~\ref{fig:LCexample}. 

Let $q$ denote the loop charge and $\dot{q}=I$ the associated loop current, then the Lagrangian for the circuit reads
\begin{equation*}
\L(q,\dot{q}) = \int \hat{\phi}(\dot{q}) \d \dot{q} - \frac{1}{2C}q^2,
\end{equation*}
which, upon substitution into (\ref{eq:EL-conservative}), yields the equation of motion $\hat\phi'(\dot{q})\ddot{q} + q/C = 0$. Note that the latter constitutes the KVL for the circuit.  

Dually, if the nonlinear current-controlled inductor is replaced by a nonlinear flux-controlled inductor, we have to consider
\begin{equation*}
\L^*(\phi,\dot{\phi}) = \frac{1}{2}C\dot{\phi}^2 - \int\hat{I}(\phi) \d \phi,
\end{equation*}
resulting in the equation of motion $C\ddot{\phi}+\hat{I}(\phi)=0$, which constitutes the KCL for the circuit.

\subsection{Conventional Non-Conservative Circuits}\label{subsec:content}

Dissipation due to conventional resistors in the circuit leads to non-conservative dynamics that is not self-adjoint. Resistors can therefore not be included using a (standard) Lagrangian function.\footnote{In some cases the dynamics can be rendered self-adjoint by looking for an integrating factor \citep{Santilli1978}. However, the form of the differential equations is altered and the Lagrangian does not have the usual interpretation of stored energy anymore \citep{Ray1979}. Furthermore, it remains unclear how to apply this method in case of complex circuits consisting of many loops.} 

To circumvent this problem, a so-called content function, which is a nonlinear multi-domain generalization of the Rayleigh dissipation function, is usually introduced.

As an example, suppose the inductor in Fig.~\ref{fig:LCexample} possesses a constant series resistance $R$. In that case, the equation of motion extends to $\phi'(\dot{q})\ddot{q} + R\dot{q} + q/C = 0$, so that $A(\dot{q})=\phi'(\dot{q})$ and $B(q,\dot{q})=R \dot{q}+q/C$. For $R = 0$, the conditions (\ref{eq:SA1})--(\ref{eq:SA4}) are clearly satisfied, but if $R \neq 0$, condition (\ref{eq:SA4}) is violated and no Lagrangian can be associated to the circuit.  On the other hand, the introduction of a function of the form $\D(\dot{q})=\frac{1}{2}R\dot{q}^2$, satisfying
\begin{equation}\label{eq:ELD}
\frac{\d}{\d t}\frac{\partial \L}{\partial \dot{q}} - \frac{\partial \L}{\partial q} + \frac{\partial \D}{\partial \dot{q}} = 0,
\end{equation}
solves the problem. However, the variational character of the dynamics that makes the Lagrangian formalism so appealing is clearly lost. 

Note that the existence of a content function relies on the fact that the resistors in the circuit are current-controlled. The inclusion of voltage-controlled resistors requires the introduction of a co-content function $\D^*(\dot\phi)$.

\section{Circuits with Memory Elements}\label{sec:EL-memory}


\subsection{Problem 1: Path-Dependence}

\begin{figure}[h]
\begin{center}
\psfrag{q}[][]{$q$}
\psfrag{p}[][]{$\phi$}
\psfrag{i}[][]{$I$}
\psfrag{C}[][]{$C$}
\psfrag{L}[r][]{$L_M$}
\includegraphics[height=30mm]{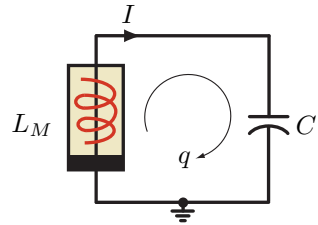}
\caption{Circuit with a meminductor.}
\label{fig:LM}
\end{center}
\end{figure}

Now consider a circuit consisting of a meminductor (\ref{eq:memL}) and a linear capacitor shown in Fig.~\ref{fig:LM}. Application of KVL yields $\dot{\phi} + q/C = 0$, or equivalently, 
\begin{equation}\label{eq:EOM-LM}
\frac{\d}{\d t}\big(L_M(q)I\big) + \frac{q}{C} = L_M(q)\dot{I} + L_M'(q)I^2 + \frac{q}{C} = 0.
\end{equation}
In order to derive the dynamics using a Lagrangian formulation, one is tempted to start from a Lagrangian that equals the magnetic co-energy stored in the meminductor minus the electric energy stored in the capacitor, i.e.,
\begin{equation}\label{eq:L-Lm}
\L(q,\dot{q}) = \frac{1}{2}L_M(q)\dot{q}^2 - \frac{1}{2C}q^2.
\end{equation}  
However, the latter is clearly not a proper state function since it depends on the path $q$.

Furthermore, substitution of the Lagrangian (\ref{eq:L-Lm}) into (\ref{eq:EL-conservative}) yields the equation $L_M(q)\ddot{q} + \frac{1}{2}L_M'(q)\dot{q}^2 + q/C = 0$, where, in comparison to (\ref{eq:EOM-LM}), we observe the appearance of an erroneous factor $\frac{1}{2}$. 

\subsection{Problem 2: Self-Adjointness}

The reason for this discrepancy is that the differential equation (\ref{eq:EOM-LM}) is not self-adjoint. To see this, let us consider the general form of a circuit consisting of memristors, meminductors, memcapacitors, and their, possibly nonlinear, conventional counterparts, independent voltage sources, and independent current sources, given by
\begin{equation}\label{eq:general-memcircuit}
A_{ij}(x,\dot{x})\ddot{x}^j + B_i(x,\dot{x})=0,   
\end{equation}
where, in comparison to (\ref{eq:general-circuit}), it is observed that $A_{ij}$ now also depends on the $x$-coordinates. The necessary and sufficient conditions for the existence of a Lagrangian for this case extend to
\begin{align}
A_{ij} &= A_{ji},\label{eq:SA1-mem}\\[0.66em]
\frac{\partial A_{ik}}{\partial \dot{x}^j} &= \frac{\partial A_{jk}}{\partial \dot{x}^i},\label{eq:SA2-mem}\\[0.3em]
\frac{\partial B_i}{\partial x^j} - \frac{\partial B_j}{\partial x^i} &= \frac{1}{2} \frac{\partial}{\partial x^k} {\left(\frac{\partial B_i}{\partial \dot{x}^j} - \frac{\partial B_j}{\partial \dot{x}^i}  \right)}\dot{x}^k ,\label{eq:SA3-mem}\\[0.3em]
\frac{\partial B_i}{\partial \dot{x}^j} + \frac{\partial B_j}{\partial \dot{x}^i} &= 2\frac{\partial A_{ij}}{\partial x^k} \dot{x}^k.\label{eq:SA4-mem}
\end{align}

Returning to the differential equation (\ref{eq:EOM-LM}), it is directly verified that, with $A(q)=L_M(q)$ and $B(q,\dot{q})=L_M'(q)\dot{q}^2+q/C$, the circuit is not self-adjoint, and therefore does not allow a Lagrangian formulation. 

A possible solution to compensate for the erroneous term is to add contra terms to the right-hand side of the differential equation; see \cite{Pesce2003}. This is tantamount to introducing a content function of the form
\begin{equation*}
\D(q,\dot{q})=\frac{1}{6}L_M'(q)\dot{q}^3,
\end{equation*}
like in (\ref{eq:ELD}). However, as argued before, the variational character is lost. 

\subsection{Conservation of Flux and Charge}

The conventional Lagrangian formalism applied to electrical circuits essentially codes the KVL and KCL in terms of energy storage via the Lagrangian. On the other hand, suppose that we integrate both Kirchhoff laws with respect to time. This would result in the laws of conservation of flux and charge around a particular loop and node
\begin{align}
\sum_i \phi_i(t) &= 0, \quad \phi_i(t) = \int\limits_{-\infty}^{t} V_i(\tau) \d \tau,\label{eq:iKVL}\\
\sum_j q_j(t) &= 0, \quad\hspace*{0.5mm} q_j(t) = \int\limits_{-\infty}^{t} I_j(\tau) \d \tau,\label{eq:iKCL}
\end{align}
respectively, which state that flux and charge can neither be created nor be destroyed \citep{Chua1969}. 

In the sequel, we will refer to (\ref{eq:iKVL}) as the iKVL and to (\ref{eq:iKCL}) as the iKCL.

\subsection{Memory State Functions}

In Section \ref{sec:memory-defs}, we have seen that the most fundamental relationship of a charge-modulated meminductor is given by $\rho=\hat{\rho}(q)$, so that (\ref{eq:memL}) is rather a consequence of the latter. It seems therefore more natural to consider, instead of the stored magnetic co-energy, which was obtained in (\ref{eq:L-Lm}) as the integral of (\ref{eq:memL}) with respect to the current, a function of the form 
\begin{equation}\label{eq:Tmem*}
\overline{T}^*(q) := \int\hat{\rho}(q)\d q.
\end{equation} 
When plotted in the $q$-versus-$\rho$ plane, (\ref{eq:Tmem*}) represents the area above the curve associated to the constitutive relation. On the other hand, its complementary part, the function $\overline{T}$, is defined to be the area below the curve and takes the form
\begin{equation}\label{eq:Tmem}
\overline{T}(\rho) := \int\hat{q}(\rho)\d \rho.
\end{equation} 
Under the assumption that $\rho=\hat{\rho}(q)$ and $q=\hat{q}(\rho)$ are invertible (i.e., one-to-one), the two state functions can be related via the Legendre transform 
\begin{equation*}
\overline{T} + \overline{T}^* = q\rho, \quad q = \frac{\d}{\d \rho}\overline{T}(\rho), \quad \rho = \frac{\d}{\d q}\overline{T}^*(q).
\end{equation*}

In a similar fashion, for a flux-modulated memcapacitor, we propose a function of the form
\begin{equation}\label{eq:Umem*}
\overline{U}^*(\phi) := \int\hat{\sigma}(\phi)\d \phi,
\end{equation} 
representing the area above the curve associated to the constitutive relation, whereas for an integrated charge-modulated memcapacitor
\begin{equation}\label{eq:Umem}
\overline{U}(\sigma) := \int\hat{\phi}(\sigma)\d \sigma,
\end{equation} 
representing the area underneath the curve. 

Furthermore, we have
\begin{equation*}
\overline{U} + \overline{U}^* = \phi\sigma, \quad \phi=\frac{\d}{\d \sigma}\overline{U}(\sigma), \quad \sigma=\frac{\d}{\d \phi}\overline{U}^*(\phi).
\end{equation*}
Note that $\overline{T}$, $\overline{T}^*$, $\overline{U}$, and $\overline{U}^*$ all serve as a state function and all exhibit the units of energy times time-squared (Joule-second-squared [Js$^2$]), which is equivalent to action times time.

\subsection{Meminductor and Memcapacitor Circuits} 

Consider again the circuit of Fig.~\ref{fig:LM}. First, we select, instead of a loop charge $q$, the integrated loop charge $\sigma$, with $\dot{\sigma}=q$, as the configuration variable. Secondly, it is observed that the conventional capacitor can be considered as a linear memcapacitor with constitutive relation $\phi=\sigma/C$. In terms of the state functions proposed in the previous subsection, let us define a Lagrangian of the form
\begin{equation*}
\overline{\L}(\sigma,\dot{\sigma}) = \int\limits_{0}^{\dot\sigma}\hat{\rho}(q)\d q - \frac{1}{2C}\sigma^2.
\end{equation*}
Then, it can be demonstrated that, invoking Hamilton's principle of least action and considering variations in terms of $\sigma$, we get the Lagrangian type of equation
\begin{equation}\label{eq:EL-sigma-example}
\frac{\d}{\d t}\frac{\partial \overline{\L}}{\partial \dot{\sigma}} - \frac{\partial \overline{\L}}{\partial \sigma} = 0,
\end{equation}  
which, in turn, generates the nonlinear differential equation $\hat{\rho}'(\dot{\sigma})\ddot{\sigma}+\sigma/C=0$. Differentiating the latter with respect to time yields
\begin{equation*}
\hat{\rho}'(\dot\sigma){\dddot{\sigma}} + \hat{\rho}''(\dot\sigma)\ddot{\sigma}^2 + \dot\sigma/C = 0, 
\end{equation*}
which, after a change of variables $\dot\sigma=q$, and recognizing that the incremental meminductance $ L_M(q) := \hat{\rho}'(\dot\sigma)$, clearly coincides with the correct equation of motion (\ref{eq:EOM-LM}). Furthermore, it is easily seen that (\ref{eq:EL-sigma-example}) constitutes the iKVL for the circuit. Thus, we have rendered the dynamics self-adjoint by considering the circuit characteristics from the perspective of flux conservation.

In general, the Lagrangian equations for circuits containing charge-modulated meminductors, integrated charge-modulated memcapacitors, and their linear conventional counterparts, take the form
\begin{equation}\label{eq:EL-sigma}
\frac{\d}{\d t}\frac{\partial \overline{\L}}{\partial \dot{\sigma}^i} - \frac{\partial \overline{\L}}{\partial \sigma^i} = 0,
\end{equation}  
with as Lagrangian $\overline{\L}(\sigma,\dot{\sigma}) = \overline{T}^*(\dot\sigma) - \overline{U}(\sigma)$, where $\overline{T}^*(\dot\sigma)$ and $\overline{U}(\sigma)$ are now representing the sums of the individual memory storage functions associated to the meminductors and the memcapacitors in the circuit, respectively.  

Naturally, the dual version of (\ref{eq:EL-sigma}) reads
\begin{equation}\label{eq:EL-rho}
\frac{\d}{\d t}\frac{\partial \overline{\L}^*}{\partial \dot{\rho}^j} - \frac{\partial \overline{\L}^*}{\partial \rho^j} = 0,
\end{equation}  
with the co-Lagrangian $\overline{\L}^*(\rho,\dot{\rho}) = \overline{U}^*(\dot\rho) - \overline{T}(\rho)$.

\subsection{Including Memristors}

Although memristors are dynamical elements they also behave as resistors. For that reason they cannot included using a Lagrangian or co-Lagrangian alone, and therefore we need to introduce, at the cost of a underlying variational principle, a second pair of state functions that play a similar role as the content and co-content functions for conventional resistors (see Subsection \ref{subsec:content}).

For a charge-modulated memristor, we propose the state function 
\begin{equation}\label{eq:Dmem}
\overline{\D}(q) := \int \hat\phi(q) \d q,
\end{equation}
which, when plotted in the $\phi$-versus-$q$ plane, represents the area underneath the constitutive relation.

For a flux-modulated memristor, we introduce the complementary function
\begin{equation}\label{eq:Dmem*}
\overline{\D}^*(\phi) := \int \hat{q}(\phi) \d \phi,
\end{equation}   
representing the area above the curve. Note that (\ref{eq:Dmem}) and (\ref{eq:Dmem*}) exhibit the units of energy times time (Joule-second [Js]), which is equivalent to action.

\begin{figure}[h]
\begin{center}
\psfrag{q1}[][]{$\sigma^1$}
\psfrag{q2}[][]{$\sigma^2$}
\psfrag{M}[][]{$R_M$}
\psfrag{R}[][]{$R$}
\psfrag{C}[][l]{$C_M$}
\psfrag{L}[r][]{$L$}
\includegraphics[height=30mm]{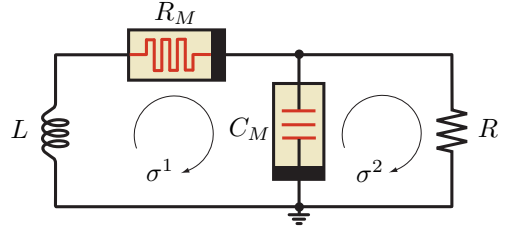}
\caption{Circuit with a memristor and a memcapacitor.}
\label{fig:RMLCm}
\end{center}
\end{figure}

As an illustration, consider the circuit depicted in Fig.~\ref{fig:RMLCm} consisting of a charge-modulated memristor $R_M$, a conventional linear inductor $L$, an integrated charge-modulated memcapacitor $C_M$, and a conventional linear resistor $R$. From a Lagrangian perspective, it is necessary to consider a loop analysis. For that, we select as configuration variables the integrated loop charges $\sigma^i$, with $i=1,2$. The Lagrangian for the circuit reads
\begin{equation*}
\overline{\L}(\sigma^1,\sigma^2,\dot{\sigma}^1) = \frac{1}{2}L(\dot\sigma^1)^2 \ - \!\!\int\limits_{0}^{\sigma^1 - \sigma^2}\hat{\phi}_{C_M}(\sigma)\d \sigma,
\end{equation*}
whereas the memristive and resistive action is given by 
\begin{equation*}
\overline{\D}(\dot\sigma^1,\dot\sigma^2) = \int\limits_0^{\dot\sigma^1} \hat\phi_{R_M}(q) \d q + \frac{1}{2}R(\dot\sigma^2)^2.
\end{equation*}

The associated Lagrangian equations take the form
\begin{align*}
\frac{\d}{\d t}\frac{\partial \overline{\L}}{\partial \dot{\sigma}^i} - \frac{\partial \overline{\L}}{\partial \sigma^i} + \frac{\partial \overline{\D}}{\partial \dot{\sigma}^i} = 0, \quad i=1,2,
\end{align*}
which provide the iKVL for loop 1:
\begin{align*}
L\ddot{\sigma}^1 + \hat{\phi}_{C_M}(\sigma^1-\sigma^2) + \hat{\phi}_{R_M}(\dot\sigma^1) = 0,
\end{align*}
and the iKVL for loop 2:
\begin{align*}
-\hat{\phi}_{C_M}(\sigma^1-\sigma^2) + R\dot\sigma^2 =0.
\end{align*}
The KVL's are obtained by differentiating the latter equations with respect to time. 

\section{Absement: a Mechanical Analogy of Integrated Charge}\label{sec:absement}

In science and engineering, the ideas and concepts developed in one branch of science and engineering are often transferred to other branches. One approach to transferring these ideas and concepts is by the use of analogies. Classically, voltage is commonly considered as the electrical analogue of force, and current is the electrical analogue of velocity. Consequently, flux can be considered as the analogue of momentum or impulse, and charge as the analogue of position or displacement.  

In mechanics, displacement and its various derivatives define an ordered hierarchy of meaningful concepts. The first derivative of displacement is velocity, the second derivative is acceleration, the third derivative is jerk, the fourth derivative is jounce, etc.. On the other hand, recently also the integral of displacement over time is introduced to model the essential behavior of flow-based musical instruments, such as the hydraulophone discussed in \cite{Mann-Absement}. This quantity is called \emph{absement}, a contraction of absence and displacement, which, in SI units, is measured in meter times seconds [ms]. One meter-second corresponds to being absent one meter from an origin or other reference point for a duration of one second. Thus, integrated charge can be considered as the electrical analogue of absement.

\subsection{Absement Related to Memcapacitors}

To gain some intuition of absement in the context of the memory elements, consider the commonly used analogy of an electrical capacitor and a bucket that can be filled with water. If current were the flow rate of water from a tap, then charge would be the amount of water in the bucket. As the flow rate is proportional to how far open the tap is, the rate of flow increases as the tap is opened up further. Now, if we consider the displacement of the tap handle from its rest position, then the amount of water in the bucket is approximately proportional to the time-integral of the handle's displacement, i.e., the handle's absement, which is a measure of how `absent' (how far and for how long) the handle is from its closed position. Equivalently, since displacement is the time-derivative of absement, the position of the handle is the time-derivative of how much water is accumulated in the bucket. 

From an electrical perspective, the absement of a capacitor is a measure for the amount of charge that is needed in a particular time interval to charge the capacitor to a certain level. For a memcapacitor this process involves both a charging and a shaping process of the capacitor and its capacitance. The amount of absement may then be used as a measure of how much charge and time is needed to write the memory of a memcapacitor. 

\subsection{Absement Related to Memristors and Meminductors}

A mechanical example that exhibits meminductive phenomena concerns the elementary problem of a heavy cable that is deployed from a reel. As argued in \cite{JeltsemaCDC2010}, under some reasonable assumptions on the geometry, this system admits a nonlinear constitutive relationship between integrated angular momentum $p$ (the mechanical analogue of integrated flux) and angular displacement (the mechanical analogue of charge) $\sigma$ of the form $\rho = \hat{\rho}(\theta)$, with $\dot \rho = p$. The amount of cable deployed from the reel equals $R\theta$, where $R$ is the diameter of the reel. The absement in this case could be a measure for the amount of cable that needs to be deployed from the system in a particular time interval to empty the reel. 

The simplest physical example of a mechanical memristor is a tapered dashpot. This type of dashpot is a mechanical resistor whose resistance (i.e., its friction coefficient) depends on the relative displacement of its terminals. Although a physical electrical passive two-terminal memristive device was constructed only recently by \cite{DmitriNature2008}, a tapered dashpot was already brought forward in the early seventies; see \cite{CSMpaper}. As for the meminductor, the amount of absement in this case could be a measure for the piston being absent from a certain resistance value for a duration of one second.

\section{Final Remarks}\label{sec:conclusions}

In this paper we have presented a novel (co-)Lagrangian framework to include memristors, meminductors, and memcapacitors, together with the conventional linear counterparts. In case of a circuit that consist only of meminductors, memcapacitors, and their conventional linear counterparts, the Lagrangian equations can be obtained from Hamilton's principle of least action. Memristive and linear resistive elements can be included via the introduction of an action function that plays a role similar to the Rayleigh (co-)dissipation function or (co-)content. We conclude the paper with the following remarks.
\begin{itemize}
\item Independent voltage (resp.~current) sources can be included in the framework by considering them either as time-varying capacitors (resp.~inductors) via a function function of the form (\ref{eq:Umem}) (resp.~(\ref{eq:Tmem})), or as time-varying resistors via a memristive action (resp.~co-action) function (\ref{eq:Dmem}) (resp.~(\ref{eq:Dmem*})). Furthermore, in order the bring the sources to the iKVL or iKCL level, we need to consider their time-integrals. 
For example, suppose we add a voltage source to the first loop of the circuit of Fig.~\ref{fig:RMLCm}. Then we need to extract either a term $-\sigma_1 \phi_e(t)$ from the Lagrangian or add a term $-\dot\sigma_1 \phi_e(t)$ to the memristive action function, with 
\begin{equation*}
\phi_e(t)=\int\limits_{-\infty}^{t} e(\tau) \d \tau,
\end{equation*} 
where $e$ denotes the source voltage. 

\item The amount of integrated charge (electrical absement) can be used as a measure of how much charge and time is needed to write the memory of a memristor, meminductor, or memcapacitor. Although not addressed in the present paper, a similar interpretation should apply to integrated flux, which in the mechanical domain corresponds to integrated momentum measured in Newton times second-squared [Ns$^2$] or kilogram-meter [kg$\cdot$m]. 

\item Although we have restricted our attention to circuits that can be modeled via either a loop or a node analysis, mixed formulations are of course also possible; see e.g., \cite{CSMpaper} for a discussion in the context of conventional circuits.
\end{itemize}

\bibliography{mathmod2012-ref}             

\end{document}